\documentclass[12pt]{article}
\usepackage{amssymb}
\usepackage{amsfonts}
\usepackage{times}
\usepackage{mathptmx}
\usepackage{amsmath}
\usepackage[usenames]{color}
\usepackage{mathrsfs}
\usepackage{amsfonts}
\usepackage{amssymb,amsmath}
\usepackage{CJK}
\usepackage{cite}
\usepackage{cases}
\usepackage{amsthm}
\def\dis{\displaystyle}\def\ssc{\scriptscriptstyle}\def\OP#1{\raisebox{-7pt}{$\stackrel{\dis\oplus}{\ssc #1}$}}

\def\cl{\centerline}

\def\a{\alpha}

\def\b{\beta}
\def\vs{\vspace*}
\def\CD{{\mathcal C}{\rm Der} }
\def\CI{{\mathcal C}{\rm Inn}}
\def\CW{\mathcal {C\!W}\!}
\def\W{\mathcal{W}}
\def\Z{\mathbb{Z}}

\def\C{\mathbb{C}}

\numberwithin{equation}{section}
\newtheorem{theo}{Theorem}[section]
\newtheorem{defi}[theo]{Definition}
\newtheorem{coro}[theo]{Corollary}
\newtheorem{lemm}[theo]{Lemma}
\newtheorem{prop}[theo]{Proposition}

\newtheorem{case}{Case}

\newtheorem{subcase}{Subcase}

\begin{document}
\begin{center}
{\bf\large A new class of  $\Z$-graded Lie conformal algebras of infinite rank}
\footnote { Supported by NSF grant Nos. 11501417, 11371278, 11431010 and  the Innovation Program of Shanghai Municipal Education Commission.

$^{\,3}$Corresponding author:  J.~Han. (jzhan@tongji.edu.cn)
}
\end{center}

\cl{Guangzhe Fan$^{\,1}$, Qiufan Chen$^{\,2}$, Jianzhi Han$^{\,1,\,3}$
}

\cl{\small $^{\,1}$School of Mathematical Sciences, Tongji University, Shanghai 200092, China}
\cl{\small $^{\,2}$Department of Mathematics, Shanghai Maritime University, Shanghai 201306, China}

\vs{8pt}
\cl{\small Emails: yzfanguangzhe@126.com, chenqf@shmtu.edu.cn,  jzhan@tongji.edu.cn}
\vs{10pt}

{\small\footnotesize
\parskip .005 truein
\baselineskip 3pt \lineskip 3pt
\noindent{{\bf Abstract:} In this paper, a new class of $\Z$-graded Lie conformal algebras $\CW(a,c)$ of  infinite rank is constructed.
 The conformal derivations and one-dimensional central extensions of $\CW(a,c)$ are completely
determined. And   all conformal modules of rank one over $\CW(a,c) (a\neq0)$ are proved to be trivial and all such nontrivial (irreducible) modules over $\CW(0,c)$ are classified.   \vs{5pt}

\parskip .001 truein\baselineskip 6pt \lineskip 6pt

}}
\parskip .001 truein\baselineskip 6pt \lineskip 6pt

\section{Introduction}
 As an algebraic language
describing the singular of the operator product expansion and a basic tool for the construction of free field realization of infinite dimensional Lie (super)algebras in conformal field theory, the notions of conformal algebras and conformal modules were initially introduced by  Kac  \cite{K3} in 1996.  Since then, the structure theory, representation theory  and cohomology theory of finite (i.e., finitely generated as $\C[\partial]$-modules)
Lie conformal algebras had been well developed (cf.  \cite{BKV, CK,CKW,DK}), and finite simple Lie conformal
superalgebras were  classified  in \cite{FK} and their representation theories were studied in \cite{BKL1,BKL2,KO}. It is natural and necessary to study the theory of infinite Lie conformal algebras. But this is definitely more challenging than that for finite case.

Important ingredients of infinite
 Lie conformal algebras are $\Z$-graded of  infinite rank. One well-known way of constructing of these Lie conformal algebras is to consider loop algebras of Lie algebras (cf. \cite{CHSX, DH,FSW,WCY}).
In this paper we construct a new class of  $\Z$-graded Lie conformal algebras  of infinite rank,      which  contains the  loop Virasoro Lie conformal algebra  \cite{WCY} and the   Block type Lie conformal algebra \cite{GXY}.   Namely, we associate  a pair of complex numbers $(a,c)$   with  a Lie conformal algebra $\CW(a,c)$ which has a $\C[\partial]$-basis
$\{L_{i}\,|\,{i}\in\Z\}$ satisfying the following $\lambda$-brackets:
\begin{eqnarray}\label{1.3}
\aligned
&[L_{i}\, {}_\lambda \, L_{j}]=\Big((ai+c)\partial+(a(i+j)+2c)\lambda\Big) L_{i+j},\quad {\rm where}\  i,j\in\Z.
\endaligned
\end{eqnarray}And this Lie conformal algebra can also be constructed by considering the Lie algebra $\W\!(a,c)$ (see Section 2).  One  motivation of the present paper is to construct new Lie conformal algebras, and another one is to  unify some known Lie conformal algebras into a class and study this class in a uniform way.

In the present paper,  the structure, representation and cohomology theories of $\CW(a,c)$ are investigated. More precisely, we show that conformal derivations of $\CW(a,c)$ are inner except for the case $a^{-1}c\in\Z$, under which case the outer conformal derivation space is one-dimensional,  and determine over which algebra there exist nontrivial conformal modules of rank one and give all nontrivial (irreducible) modules of such kind. Finally, we compute one-dimensional central extensions of $\CW(a,c)$ and  determine all such extensions. In particular, some of the results in  \cite{GXY,WCY} are recovered here.

The following  is the organization of this paper. In Section $2$, some basic definitions of Lie conformal algebras are reviewed. Moreover, we construct a new
class of $\Z$-graded Lie conformal algebras of infinite rank, denoted by $\CW(a,c)$. In Section $3$, we determine the conformal derivations of $\CW(a,c)$. Section 4 is devoted to giving a classification of nontrivial (irreducible) conformal modules of rank one.   Finally, we give all one-dimensional
central extensions of $\CW(a,c)$ in Section 5.

Throughout the paper, we denote by $\C,\,\C^*,\, \Z$ the sets of all complex numbers, nonzero complex numbers, integers,
respectively, and all vector spaces, tensor
products are considered over $\C$ and all linear maps are $\C$-linear.

\section{Preliminaries and constructions of $\CW(a,c)$}
In this section, we first recall some definitions related to Lie conformal algebras  (cf. \cite{K1,K3}) and then construction a Lie conformal algebra $\CW(a,c)$ for any pair $(a,c)$ of complex numbers.

\begin{defi}\label{210}\rm
A Lie conformal algebra is a $\C[\partial]$-module $R$ endowed with a linear map $R\otimes R\rightarrow R[\lambda]$, $a\otimes b\rightarrow [a{}\,
_\lambda \, b]$, called $\lambda$-bracket, where $\lambda$ is an indeterminate variable and $R[\lambda]=\C[\lambda]\otimes R$, subject to the following
three axioms:
\begin{eqnarray}\label{211}
&& {\rm(conformal\ sesquilinearity)}\quad [\partial a\,{}_\lambda \,b]=-\lambda[a\,{}_\lambda\, b],
[a\,{}_\lambda \,\partial b]=(\partial+\lambda)[a\,{}_\lambda\, b];\\
&& {\rm (skew\ symmetry)}\quad\quad \quad \quad~~~~ [a\, {}_\lambda\, b]=-[b\,{}_{-\lambda-\partial}\,a];\\
&&{\rm (Jacobi\ identity)}\quad\quad\quad \quad ~~~~~~[a\,{}_\lambda\,[b\,{}_\mu\, c]]=[[a\,{}_\lambda\, b]\,{}_{\lambda+\mu}\, c]+[b\,{}_\mu\,[a\,{}_\lambda \,c]].
\end{eqnarray}
\end{defi}

\begin{defi}\label{212}\rm
A conformal module over a Lie conformal algebra $R$ or a $R$-module  is a $\C[\partial]$-module $V$ endowed with a $\lambda$-action $R\otimes V\rightarrow V[\lambda]$ such
that for any $a,b\in R$ and $v\in V$:
\begin{equation*}\label{213}
\aligned
&(\partial a)\,{}_\lambda\, v=-\lambda a\,{}_\lambda\, v,\ \ \ \ \ a{}\,{}_\lambda\, (\partial v)=(\partial+\lambda)a\,{}_\lambda\, v,\\
&a\,{}_\lambda\, (b{}\,_\mu\, v)-b\,{}_\mu\,(a\,{}_\lambda\, v)=[a\,{}_\lambda\, b]\,{}_{\lambda+\mu}\, v.
\endaligned
\end{equation*}
If, in addition,  $V$ is a free $\C[\partial]$-module of rank $n$, then we say $V$  a $R$-module of  rank $n$.
\end{defi}
\begin{defi}\label{214}\rm
A Lie conformal algebra $R$ is {\it $\Z$-graded} if $R=\oplus_{i\in \Z}R_i$, where each $R_i$ is a $\C[\partial]$-submodule
such that $[R_i\,{}_\lambda\, R_j]\subset R_{i+j}[\lambda]$ for any $i,j\in \Z$.
Similarly, a $R$-module $V$  is  {\it $\Z$-graded} if $V=\oplus_{i\in \Z}V_i$, where each $V_i$ is a $\C[\partial]$-submodule satisfying
$R_i\,{}_\lambda\, V_j\subset V_{i+j}[\lambda]$ for any $j\in \Z$.
\end{defi}

\begin{defi}\label{215}\rm
Let $V$ and $W$ be two $\C[\partial]$-modules. A conformal linear map from $V$ to $W$ is a linear map $\phi_\lambda:V\longrightarrow
\C[\partial][\lambda]\otimes_{\C[\partial]}W$ such that
\begin{equation*}
\phi_\lambda(\partial v)=(\partial+\lambda)\phi_\lambda(v)\ \ \mbox{ \ for $v\in V$.}
\end{equation*}
\end{defi}

In the rest part of this section we are going to give  constructions of  $\Z$-graded Lie conformal algebras $\CW(a,c)$ and Lie algebras $\W\!(a,c)$.

Consider first the $\C[\partial]$-module $\CW$ which has  a $\C[\partial]$-basis $\{L_{i}\,|\,{i}\in\Z\}.$ The goal of this section is to make $\CW$ carry the structure of a Lie conformal algebra such that the $\lambda$-brackets have the following form:
\begin{eqnarray}\label{3.1}
\aligned
&[L_{i}\, {}_\lambda \, L_{j}]=(f(i,j)\partial+g(i,j)\lambda) L_{i+j}\quad {\rm for}\ i,j\in\Z,
\endaligned
\end{eqnarray}
where $f(x,y),g(x,y)\in\C[x,y]$.

The conformal skew symmetry requires
\begin{equation*}
[L_{j}\, {}_\lambda \, L_{i}]=-[L_{i}\,{}_{-\lambda-\partial}\,L_{j}],
\end{equation*}
i.e.,
\begin{equation*}
f(j,i)\partial+g(j,i)\lambda =-\Big(f(i,j)\partial+g(i,j)(-\lambda-\partial)\Big)
\end{equation*}hold for any $i,j\in\Z$. This is equivalent to
\begin{equation}\label{3.4}
\aligned
&g(i,j) =f(i,j)+f(j,i)\quad {\rm for}\ i,j\in\Z.
\endaligned
\end{equation}Taking this into consideration,
(\ref{3.1}) can be written as
\begin{eqnarray}\label{3.6}
\aligned
&[L_{i}\, {}_\lambda \, L_{j}]=\Big(f(i,j)\partial+(f(i,j)+f(j,i))\lambda\Big) L_{i+j}.
\endaligned
\end{eqnarray}
Now  the conformal Jacobi identity condition
$$[L_{i}\,{}_\lambda\,[L_{j}\,{}_\mu\, L_{k}]]=[[L_{i}\,{}_\lambda\, L_{j}]\,{}_{\lambda+\mu}\, L_{k}]+[L_{j}\,{}_\mu\,[L_{i}\,{}_\lambda
\,L_{k}]$$
is equivalent to
\begin{equation}\label{3.8}
\aligned
&f(i,k)f(j,i+k) =f(j,k)f(i,j+k)
\endaligned
\end{equation}
and
\begin{equation}\label{3.9}
\aligned
&f(j,i)(f(i+j,k)+f(k,i+j)) =f(j,k)(f(i,j+k)+f(j+k,i)).
\endaligned
\end{equation}

\begin{lemm}\label{theo3.1}
Let $f(x,y), g(x,y)\in\C[x,y]$ be as above.  Then $$f(x,y)=ax+c\ \ and\ \ g(x,y)=a(x+y)+2c\quad { for\ some}\ a,c\in\C.$$
\end{lemm}
\begin{proof} Note from \eqref{3.8} that ${\rm deg}_{x}f(x,y)+{\rm deg}_{y}f(x,y)={\rm deg}_{x}f(x,y)$ by comparing the degree of $i$.  Thus,   $f(x,y)\in\C[x].$ It follows from this and setting $k=0$ in \eqref{3.9},
and one can see that ${\rm deg}_{x}f(x,y)\le 1$. So  $f(x,y)$ and thereby $g(x,y)$ by \eqref{3.4} have the promised forms.
\end{proof}

Now we arrive at the following result.

\begin{prop}\label{theo3.2}
Let $\CW$ be  a free $\C[\partial]$-module with basis $\{L_{i}\,|\,{i}\in\Z\}$. Then $\CW$ carries the structure of a Lie conformal algebra whose $\lambda$-bracket has the form
\begin{equation*}
[L_{i}\, {}_\lambda \, L_{j}]=(f(i,j)\partial+g(i,j)\lambda) L_{i+j} \ for\ some f(x,y),g(x,y)\in\C[x,y],
\end{equation*} if and only if there exist $a,c\in\C$ such that
\begin{eqnarray}\label{3.15}
&[L_{i}\, {}_\lambda \, L_{j}] =\Big((a i+c)\partial+(a(i+j)+2c)\lambda\Big) L_{i+j}.
\end{eqnarray} \end{prop}

Denote  the Lie conformal algebra  in Proposition \ref{theo3.2}  by $\CW(a,c)$, which is $\Z$-graded (cf. Definition \ref{214}):  $\CW(a,c)=\oplus_{{i}\in\Z} \ \CW(a,c)_{i}$ with
$\CW(a,c)_{i}=\C[\partial]{L_{i}}$. Thus, we obtain a class of $\Z$-graded Lie conformal algebras of infinite rank. In particular,   the formula (\ref{3.15}) becomes
\begin{eqnarray}\label{3.17}
\aligned
&[L_{i}\, {}_\lambda \, L_{j}] =(\partial+2\lambda) L_{i+j},
\endaligned
\end{eqnarray} when $a=0, c=1$ and $\CW(0,1)$ turns out to be  the loop Virasoro Lie conformal algebra studied in \cite{WCY}; while $a=1,c=0$ the formula (\ref{3.15}) is simplified as
\begin{eqnarray}\label{3.16}
\aligned
&[L_{i}\, {}_\lambda \, L_{j}] =\Big(i\partial+(i+j)\lambda\Big) L_{i+j},
\endaligned
\end{eqnarray} and  $\CW(1,0)$ is the Lie conformal algebra of a Block type Lie algebra (cf.\cite{GXY}). However, the Lie conformal algebras $\CW(a,c)$ for many other pairs $(a,c)$ are first introduced in the present paper, that is the reason that  this class of Lie conformal algebra is called a new one.

Another construction of $\CW(a,c)$  is given from the view point of Lie algebras.
Let $\W$ be the vector space with basis $\{L_{i,p}\mid i,p\in\Z\}$. For any ${i}\in\Z$, set $$L_{i}(z)=\sum_{p\in\Z} L_{i,p}z^{-p-2}.$$ For any $i,j\in\Z$, define \begin{eqnarray*}&&
[L_{i}(z),L_{j}(w)]=(ai+c)\partial_w L_{i+j}(w)\delta(z,w)+(a(i+j)+2c)L_{i+j}(w)\partial_w\delta(z,w),
\end{eqnarray*} which is equivalent to \begin{eqnarray}\label{eq-add-1}
[L_{i,p},L_{j,q}]=\Big(a(j(p+1)-i(q+1))+c(p-q)\Big)L_{i+j,p+q}\quad{\rm for}\ p,q\in\Z.
\end{eqnarray} Then under the Lie brackets given by \eqref{eq-add-1}, one can check that $\W$ carries the structure of a Lie algebra which is denoted by $\W\!(a,c)$. Note that the Lie algebra $\W\!(a,c)$ and the Lie conformal algebra $\CW(a,c)$ are associated with each other. So this association gives another construction of $\CW(a,c).$

\section{Conformal derivations of $\CW(a,c)$}
\begin{defi}\label{40}\rm
Let $R$ be a Lie conformal algebra. A conformal linear map $D_\lambda:R\longrightarrow R[\lambda]$ is called a conformal derivation if
\begin{equation*}\label{41}
D_\lambda([a\,{}_\mu \,b])=[(D_\lambda a)\,{}_{\lambda+\mu} \,b]+[a\,{}_\mu \,(D_\lambda b)]\ \ \mbox{ \ for $a,b\in R$.}
\end{equation*}If in addition, $D_\lambda=({\rm adx}_x)_\lambda$ for some $x\in R$ with $({\rm ad}_x)_\lambda y= [x\, {}_\lambda\, y]$ for $y\in R,$ then $D_\lambda$ is called a conformal inner derivation.
\end{defi}
Denote by $\CD(\CW(a,c))$ and $\CI(\CW(a,c))$ the vector spaces of all conformal derivations and conformal inner  derivations of $\CW(a,c)$, respectively. We write a conformal derivation $D$ rather than $D_{\lambda}$ for convenience.  For any $D\in
 \CD(\CW)$, define $D^j(L_i)=\pi_{i+j} D(L_i)$ for any $i\in\Z$, where
$$\pi_{j}: \C[\lambda]\otimes \CW\cong \OP{k\in\Z}\C[\partial,\lambda]L_k \rightarrow \C[\partial,\lambda]{L_{j}}$$ is the natural projection.
Then $D^j$ is a conformal derivation and $D=\sum_{j\in\Z} D^j$ in the sense that for any $x\in \CW$ only finitely many of $D^j_\lambda(x)$ for $j\in\Z$ are nonzero.

\begin{theo}\label{main-5}
We have \begin{eqnarray*}\CD(\CW(a,c))=\begin{cases}\CI(\CW(a,c)), &\mbox{if}\   a=0\ {\rm or}\ a^{-1}c\notin\Z,\\
\CI(\CW(a,c))\oplus \C \mathcal D^{\frac{c}{a}}, &\mbox{if}\   a^{-1}c\in\Z,\\
\end{cases}
\end{eqnarray*}
\end{theo} where $\mathcal D^0(L_i)=L_i$ and $\mathcal D^{\frac{c}{a}}(L_i)=(1-\frac a c)L_{i+ \frac{c}{a}}$ for $i\in\Z$.
\begin{proof}
It remains to consider the case $(a,c)\neq (0,0)$, since  $\CD(\CW(0,0))=0$  by the definition of conformal derivations.

Let $D\in \CD(\CW(a,c))$ and assume that $$D^j_\lambda(L_i)=f_i^j(\partial,\lambda)L_{i+j},\quad{\rm  where}\ f^j_i(\partial,\lambda)\in\C[\partial,\lambda].$$
Applying $D^j_\lambda$ to $[L_0\ {}_\mu \ L_i]=\big(c\partial+(ai+2c)\mu \big ) L_i$, one has
\begin{eqnarray}\label{5.14}
\big(c(\partial+\lambda)+(ai+2c)\mu\big)f^j_i(\partial,\lambda)\!\!\!\!&=&\!\!\!\! f^j_0(-\lambda-\mu,\lambda)\Big((aj+c)\partial+(\lambda+\mu)(a(j+i)+2c)\Big)+\nonumber\\ &&f^j_i(\partial+\mu,\lambda)\Big(c\partial+(a(i+j)+2c)\mu\Big) \quad{\rm for}\ i,j\in\Z.
\end{eqnarray}
Setting $\mu=0$ in (\ref{5.14}), we have
\begin{equation}\label{5.15}
c\lambda f^j_i(\partial,\lambda)=f^j_0(-\lambda,\lambda)\Big((aj+c)\partial+(a(j+i)+2c)\lambda\Big) \quad{\rm for}\ i,j\in\Z.
\end{equation}

 Consider first that $c=0$. Then by taking $j=0$, \eqref{5.14} reduces to  \begin{eqnarray}\label{add-5.14}
i\big((\lambda+\mu)f^0_0(-\lambda-\mu,\lambda)-\mu f^0_i(\partial+\mu,\lambda)
-\mu f^0_i(\partial,\lambda)\big)=0\quad{\rm for}\ i\in\Z,
\end{eqnarray} from which together with \eqref{5.15} it is not hard to see that $f^0_i(\partial,\lambda)=i (\lambda g(\lambda)+b)$ for  some $g(\lambda)\in\C[\lambda], b\in\C$ and any $i\in\Z$. Then  one can check that $D^0={\rm ad}_{\frac1 a g(-\partial)L_0}+b\mathcal D^0$. While for $j\neq 0$, letting $i=0$ in \eqref{5.14} gives rise to  \begin{equation}\label{5.11}
f^j_0(-\lambda-\mu,\lambda)(\partial+\lambda+\mu)+f^j_0(\partial+\mu,\lambda)\mu=0.
\end{equation}
It follows from by respectively comparing the coefficients of $\partial$ and $\lambda$ that $f^j_0(\partial,\lambda)=0$. Whence \eqref{5.14} becomes
\begin{equation*}\label{5.12}
(i+j)f^j_i(\partial+\mu,\lambda)=i f^j_i(\partial,\lambda)\quad{\rm for}\ i,j\in\Z,
\end{equation*}
which forces  $f^j_i(\partial,\lambda)=0$. That is, $D^j_\lambda(L_i)=0$ for all $j\neq 0$ and thereby $$D=D^0={\rm ad}_{\frac1 a g(-\partial)L_0}+b\mathcal D^0\in \CI(\CW(a,c))\oplus \C \mathcal D^{0}.$$

Next we assume that $c\neq 0$. Let $j\in\Z$ be such that $aj+c\neq 0$.  Then by \eqref{5.15},  $\lambda | f^j_0(-\lambda,\lambda).$  Setting $h_j(\lambda)=\frac{f^j_0(-\lambda,\lambda)}{\lambda}$, one has $D^j={\rm ad}_{h_j(-\partial)L_j}.$ Consider that $a^{-1}c\in\Z$.  Setting $j=a^{-1}c$ in  (\ref{5.15}) we have
\begin{equation*}\label{5.16}
f^{a^{-1}c}_i(\partial,\lambda)=f^{a^{-1}c}_0(-\lambda,\lambda)(1-\frac a ci) \quad{\rm for}\ i\in\Z.
\end{equation*} Write $f^{a^{-1}c}_0(-\lambda,\lambda)=\lambda l(\lambda)+d$, where $l(\lambda)\in\C[\lambda]$ and $d\in\C$. Then $D^{\frac c a}={\rm ad}_{a^{-1}l(-\partial)L_{\frac c a}}+d\mathcal D^{\frac c a}.$

Note that   $\sum_{j\in\Z\setminus \{a^{-1}c\}}D^j=\sum_{j\in\Z\setminus \{a^{-1}c\}}{\rm ad}_{h_j(-\partial)L_j}$ is a finite sum, since otherwise a contradiction will arise: $$D_\lambda(L_0)\in\prod_{i\in\Z}\CW(a,c)_i\setminus \bigoplus_{i\in\Z}\CW(a,c)_i.$$  Then, $\sum_{j\in\Z\setminus \{a^{-1}c\}}D^j\in \CI(\CW(a,c))$,  which implies $$D=\sum_{j\in\Z}D^j\in \CI(\CW(a,c))\quad  {\rm if}\ a=0\  {\rm or}\  a^{-1}c\notin\Z$$ and $$D=\sum_{j\in\Z}D^j\in \CI(\CW(a,c))\oplus \C\mathcal D^{\frac c a}\quad {\rm if}\ a^{-1}c\in\Z.$$
\end{proof}

\section{$\CW(a,c)$-modules of rank one}

In this section we shall determine over which Lie conformal algebra $\CW(a,c)$ there exist nontrivial conformal modules of rank one and  give a classification of such nontrivial modules. The main result of this section is formulated as follows:

\begin{theo}\label{theo6.2}
Suppose that $M=\C[\partial]v$ is  a nontrivial $\CW(a,c)$-module of rank one. Then $a=0$ and the  $\lambda$-actions on $M$ are given by
\begin{equation*}
L_i\, {}_\lambda \,v=h_i(\lambda)v\quad {\rm for\ some}\ h_i(\lambda)\in \C[\lambda]\ {\rm if}\ c=0
\end{equation*} and
\begin{equation*}
L_i\, {}_\lambda \,v={\tau}^i(c\partial+\rho\lambda+\sigma)v\quad {\rm for\ some}\  \rho, \sigma, \tau\in\C \ {\rm if}\ c\neq0.
\end{equation*}
\end{theo}

\begin{proof}
Assume that $L_i\,{}_\lambda\, v=f_i(\partial,\lambda)v$ for some $f_i(\partial,\lambda)\in\C[\partial,\lambda]$ and $i\in\Z$. It follows from $[L_i\,{}_\lambda\, L_j]\,{}_{\lambda+\mu} v=L_i\,{}_\lambda\,(L_j\, {}_\mu v)+L_j\,{}_\mu\, (L_i \,{}_\lambda v)$ that

\begin{eqnarray}\label{6.6}
&&\Big((aj+c)\lambda-(ai+c)\mu)\Big)f_{i+j}(\partial,\lambda+\mu)\nonumber\\ &=&f_j(\partial+\lambda,\mu)f_i(\partial,\lambda)-f_i(\partial+\mu,\lambda)f_j(\partial,\mu)\quad {\rm for}\ i,j\in\Z.
\end{eqnarray}
Note from \eqref{6.6} that ${\rm deg}_\partial f_i(\partial,\lambda)=0$ and therefore $f_i(\partial,\lambda)=h_i(\lambda)$  for any $i\in\Z$ if $(a,c)=(0,0).$  Similar arguments as in \cite{WCY} can be applied to  the case $a=0,c\neq0$ and show  the desired $\lambda$-actions as in the statement. So it suffices to show $f_i(\partial,\lambda)=0$  for all $i\in\Z$  in the case $a\neq 0$.

Setting $i=j=0$ in \eqref{6.6} gives
\begin{equation}\label{6.8}
c(\lambda-\mu)f_0(\partial,\lambda+\mu)=f_0(\partial+\lambda,\mu)f_0(\partial,\lambda)-f_0(\partial+\mu,\lambda)f_0(\partial,\mu).
\end{equation}
Assume that $f_0(\partial,\lambda)=\sum_{k=0}^t\partial^k h_k(\lambda)$ with $h_k(\lambda)\in\C[\lambda]$ and $h_t(\lambda)\neq 0$.  Submitting this into \eqref{6.8} we see that $t\le1$. Moreover, $f_0(\partial,\lambda)=h_0(\lambda)$ if $c=0$ and
either $f_0(\partial,\lambda)=0$ or $f_{0}(\partial,\lambda)=c\partial+\rho\lambda+\sigma$ for some $\rho, \sigma\in\C$ if $c\neq0$. Taking $i=0$ in (\ref{6.6}),  we have

\begin{eqnarray}\label{6.9}
&&f_j(\partial+\lambda,\mu)(c\partial+d\rho+\sigma)-f_j(\partial,\mu)(c\partial+c\mu+\rho\lambda+\sigma)\nonumber\\ &=&((aj+c)\lambda-c\mu)f_{j}(\partial,\lambda+\mu)\quad  {\rm for}\ j\in\Z.
\end{eqnarray}
 Then comparing the degrees on $\partial$ of both sides of \eqref{6.9} gives $f_i(\partial,\lambda)=0$ for $i\in\Z$ provided $f_0(\partial,\lambda)\in\C[\lambda]$. Thus, this completes the case $c=0$ and also the case $c\neq0$ but $f_0(\partial,\lambda)=0$.

The remaining case is $a,c\in\C^*$ and $f_{0}(\partial,\lambda)=c\partial+\rho\lambda+\sigma.$ Let $j\in\Z$ be such that $aj+c\neq0$. Taking $\mu=c^{-1}(aj+c)\lambda$ in (\ref{6.9}), we obtain
\begin{equation*}\label{6.10}
f_j(\partial+\lambda,c^{-1}(aj+c)\lambda)(c\partial+\rho\lambda+\sigma)=f_j(\partial,c^{-1}(aj+c)\lambda)(c\partial+(aj+c)\lambda+\rho\lambda+\sigma),
\end{equation*}
from which we can easily see that  $f_j(\partial,\lambda)=0$ by comparing the degrees of $\partial$ and $\lambda$, respectively. That is, $f_j(\partial,\lambda)=0$ for all $j\in\Z$ such that $aj+c\neq0$. But it follows from this and \eqref{6.6} that $f_j(\partial,\lambda)=0$ also holds  for $j$ such that $aj+c=0$. This completes the last case.
\end{proof}

Denote the modules in Theorem \ref{theo6.2} by $M_{h_i(\lambda)}$ and $M_{\tau, \rho,\sigma}$, respectively.

\begin{coro}
Let $M$ be a nontrivial irreducible  $\CW(0,c)$-module of rank one. Then $M$ is isomorphic to $M_{h_i(\lambda)}$ for some $0\neq\big(h_i(\lambda)\big)_i\in\prod_i\C[\lambda]$ if $c=0$ and to $M_{\tau, \rho,\sigma}$ for some $\tau,\rho\in\C^*, \sigma\in\C$ if $c\neq 0$.
\end{coro}
\begin{proof}
 The irreducibilities of $M_{h_i(\lambda)}$ and $M_{\tau, \rho,\sigma}$ can  be checked easily. Then it follows immediately from Theorem \ref{theo6.2}.
\end{proof}

\section{One-dimensional central extensions of $\CW(a,c)$}

An extension of a Lie conformal algebra $R$ by an abelian Lie conformal algebra $G$ is a short exact sequence of Lie conformal algebras
\begin{equation*}
0\rightarrow G \rightarrow \widehat{R} \rightarrow R \rightarrow 0,
\end{equation*}
and  $\widehat{R}$ is called an extension of $R$ by $G$. This extension is said to be central if
\begin{equation*}
G\subset Z(\widehat{R})=\{x\in\widehat{R}\,|\, [x\,_\lambda \, y]=0, \forall~y\in \widehat{R}\},~\partial G=0.
\end{equation*}
Let $\widehat{R}$ be a central extension of $R$ by a one-dimensional center $\C\mathfrak{c}$, i.e.,  $\widehat{R}\cong R\bigoplus \C\mathfrak{c}$
as vector spaces and
$$[x\,_\lambda\, y]_{\widehat{R}}=[x\,_\lambda\, y]_{R}+\phi_{\lambda}(x, y)\mathfrak{c}$$
for some {\em $2$-cocycle $\phi_\lambda$ of $R$}, by which we mean that $\phi_\lambda: R\otimes R \rightarrow \C[\lambda]$ is a bilinear map satisfying
\begin{eqnarray}\label{701}
&&({\rm skew\text{-}symmetry})\quad\quad\quad\quad~~~~ \phi_{\lambda}(x, y)=- \phi_{-\lambda}(y, x)\nonumber ,\\
&& {\rm (conformal\ sesquilinearity)}\quad\phi_{\lambda}(\partial x, y)=-\lambda \phi_{\lambda}(x, y)=-\phi_{\lambda}(x,\partial y),\\
&&{\rm (Jacobi\ identity)}\quad\quad\quad\quad\quad~\phi_{\lambda+\mu}([x\,_\lambda \, y], z)=\phi_{\lambda}(x, [y\,_\mu\, z])-\phi_{\mu}(y, [x\,_\lambda\, z])\quad {\rm for}\ x,y,z\in R.\nonumber
\end{eqnarray}

This section is devoting to determining all the central extensions $\widehat{\CW(a,c)}$ of $\CW(a,c)$ by a one-dimensional center $\C\mathfrak{c}$, i.e., $\widehat{R(a,c)}=R(a,c)\oplus\C\mathfrak{c}$, on which the $\lambda$-brackets (\ref{1.3}) are replaced by
\begin{eqnarray*}&&
[L_i\,{}_\lambda\, L_j]=\Big((ai+c)\partial+(a(i+j)+2c)\lambda\Big)L_{i+j}+\phi_{\lambda}(L_{i},L_{j})\mathfrak{c}
\end{eqnarray*}
for some 2-cocycle $\phi_\lambda$.  Thus,  it is sufficient to determine all 2-cocycles of $\CW(a,c)$. Observe from the conformal sesquilinearity for $\phi_\lambda$ that   the values  $\phi_\lambda(x,y)$ for arbitrary $x,y\in \CW(a,c)$ are uniquely determined by $\phi_\lambda(L_i,L_j)$ for all $i,j\in\Z$. Thus it is sufficient to compute $\phi_\lambda(L_i,L_j)$. Consider first the case $(a,c)=(0,0)$. Then  the Jacobi identity will automatically hold and a 2-cocycle  is a bilinear map satisfying the skew-symmetry. So we only focus on the case $(a,c)\neq (0,0)$.

\begin{theo}\label{7000}
Let $\widehat{\CW(a,c)}$ be a  one-dimensional central extension  of $\CW(a,c).$ Then   the $\lambda$-bracket on $\widehat{\CW(a,c)}$ for $i,j\in\Z$ has the following form:
\begin{eqnarray*}&&
[L_i\,{}_\lambda\, L_j]=c(\partial+2\lambda)L_{i+j}+(A(i+j)\lambda+B(i+j)\lambda^3)\mathfrak{c}
\end{eqnarray*}for some complex functions $A, B$ (i.e., $A$ is  a map from $\Z$ to $\C$) if $a=0$;
\begin{eqnarray*}
[L_i\,{}_\lambda\, L_j]=\Big((ai+c)\partial+(a(i+j)+2c)\lambda\Big) L_{i+j}+(i\delta_{c,0}\delta_{i+j,0}\Delta+A(i+j)\lambda)\mathfrak{c}
\end{eqnarray*}for some complex function $A$, $\Delta\in \C$ if $a\in\C^*$.

\end{theo}

\begin{proof}
Let $\phi_\lambda$ be a 2-cocycle of $\CW(a,c).$    Assume that $\phi_{\lambda}(L_{i},L_{j})=\sum_{m}a_{m}(L_{i},L_{j})\lambda^{m}\in \C[\lambda]$. It follows from applying the Jacobi identity  in \eqref{701} to the triple $(L_{i},L_{j},L_{k})$ that
\begin{eqnarray*}\label{711}
\Big((aj+c)\lambda-(ai+c)\mu\Big)\phi_{\lambda+\mu}(L_{i+j}, L_{k})\!\!\!&=&\!\!\!\Big((aj+c)\lambda+(a(j+k)+2c)\mu\Big)\phi_{\lambda}(L_{i}, L_{j+k})-\nonumber\\
&&\Big((ai+c)\mu+(a(i+k)+2c)\lambda\Big)\phi_{\mu}(L_{j}, L_{i+k})
\end{eqnarray*}
and thereby
\begin{eqnarray}\label{712}
&&\Big((aj+c)\lambda-(ai+c)\mu\Big)\sum_{m}a_{m}(L_{i+j},L_{k})(\lambda+\mu)^{m}\nonumber\\ &=&
\Big((aj+c)\lambda+(a(j+k)+2c)\mu\Big)\sum_{m}a_{m}(L_{i},L_{j+k})\lambda^{m}-\\ &&\Big((ai+c)\mu+(a(i+k)+2c)\lambda\Big)\sum_{m}a_{m}(L_{j},L_{i+k})\mu^{m}\nonumber.
\end{eqnarray}
Note from the latter formula that ${\rm deg}_\lambda \phi_\lambda(L_i,L_j)\le3$  for any $i,j\in\Z$ and the equality only occurs when $a=0$. Whence the discussion is divided into the following two cases. In what follows,  $L_k$ is treated as zero if $k\notin\Z$.

\begin{case}\label{caseee-1}
$a\neq0.$
\end{case}

Note that ${\rm deg}_\lambda\phi_\lambda(L_i,L_j)\le 2$ in this case. We first consider the values $a_2(L_i, L_j)$ for $i,j\in\Z$. Comparing the coefficients of $\lambda^3$ on both sides of \eqref{712} we have $a_2(L_{i+j}, L_k)=a_2(L_i, L_{j+k})$ unless $j=-a^{-1}c.$   From this, it is not hard to derive $a_2(L_{i+j}, L_k)=a_2(L_i, L_{j+k})$ for all $i,j,k\in\Z$, which is equivalent to saying the values $a_2(L_i, L_j)$ depend only on the sum $i+j$. But now taking the skew-symmetry for $\phi_\lambda$ into account  forces $a_2(L_i, L_j)=0$ for all $i,j\in\Z$. Similarly, we can show that the values $a_1(L_i, L_j)$ depend only the values $i+j$, which allows us to define a complex function $A: \Z\rightarrow \C$ given by sending $i+j$ to $a_1(L_i,L_j).$
It follows from comparing the coefficients of $\lambda$ in  \eqref{712} that \begin{equation}\label{713-1}
(aj+c) a_{0}(L_{i+j},L_{k})=(aj+c) a_{0}(L_{i},L_{j+k})-\big(a(i+k)+2c\big)a_{0}(L_{j},L_{i+k}).
\end{equation} In particular, taking $j=0$ in  \eqref{713-1} gives rise to
$a_{0}(L_{0},L_{k})=0$ for all $k\neq -2a^{-1}c.$

\begin{subcase}
$c\neq0.$
\end{subcase}
Note from setting $k=-2a^{-1}c-i$ in \eqref{713-1} that $a_0(L_{i+j}, L_{-2a^{-1}c-i})=a_0(L_i,L_{-2a^{-1}c-i+j})$ for any $j\neq-a^{-1}c$. In fact, one can deduce that the above formula also holds for $j=-a^{-1}c$. Thus, the values $a_0(L_i,L_j)$ for $i+j=-2a^{-1}c$ are equal to a constant, say,   $\xi$.
Setting $i=0$ in \eqref{713-1} gives
\begin{equation}\label{714---2}
\big(a(j+k)+3c\big)a_0(L_j,L_k)=(aj+c)a_0(L_0, L_{j+k})\quad{\rm for}\ j,k\in\Z.
\end{equation}
In particular, letting $k=-2a^{-1}c-j$ in \eqref{714---2} we have

\begin{equation*}\label{714}
c\xi=(aj+c)\xi,
\end{equation*} which implies $\xi=0$, namely, $a_0(L_i, L_j)=0$ for $i+j=-2a^{-1}c$.  Now   taking $j+k=-2a^{-1}c$ in \eqref{713-1} we see that $a_0(L_j,L_k)=0$ unless $j+k=-3a^{-1}c$. Similarly, it can be shown that $a_0(L_j,L_k)=0$ for all $j,k\in\Z$.
\begin{subcase}
$c=0.$
\end{subcase}
In this case, $a_0(L_0,L_k)=0$ for all $k\in\Z$ by using the skew-symmetry and then $a_0(L_i,L_j)=0$  unless $i+j=0$.
Setting $k=-i-j$ in (\ref{713-1}), one has
\begin{equation*}\label{715}
j a_{0}(L_{i+j},L_{-i-j})=j a_{0}(L_{i},L_{-i})+j a_{0}(L_{j},L_{-j}),
\end{equation*}
from which and the skew-symmetry  we obtain
\begin{equation*}\label{716}
a_{0}(L_{i},L_{-i})=i a_{0}(L_{1},L_{-1})\quad \mbox{ \ for } i\in\Z.
\end{equation*}
To sum up,  $a_{0}(L_{i},L_{j})=i\delta_{i+j,0} \Delta$ for $i\in\Z$,  where $\Delta=a_{0}(L_{1},L_{-1})$.

\begin{case}
$a=0.$
\end{case}

In this case, \eqref{712} turns out to be \begin{eqnarray}\label{712-1}
&&(\lambda-\mu)\sum_{m=0}^3a_{m}(L_{i+j},L_{k})(\lambda+\mu)^{m}\nonumber\\ &=&
(\lambda+2\mu)\sum_{m=0}^3a_{m}(L_{i},L_{j+k})\lambda^{m}-(\mu+2\lambda)\sum_{m=0}^3a_{m}(L_{j},L_{i+k})\mu^{m}.
\end{eqnarray}
Comparing the coefficients of $\lambda^4$ gives  $a_3(L_{i+j}, L_\gamma)=a_3(L_i,L_{j+k})$ for any $i,j,k\in\Z$.  And this entails us to define a complex function $B: \Z\rightarrow\C$ sending $i+j$ to $a_3(L_i,L_j),$ i.e., $B(\a+\b)=a_3(L_i,L_j)$ for $i,j\in\Z$. Note  that \eqref{712-1} can be reduced to \begin{eqnarray}\label{712-2}
&&(\lambda-\mu)\sum_{m=0}^2a_{m}(L_{i+j},L_{k})(\lambda+\mu)^{m}\nonumber\\ &=&
(\lambda+2\mu)\sum_{m=0}^2a_{m}(L_{i},L_{j+k})\lambda^{m}-(\mu+2\lambda)\sum_{m=0}^2a_{m}(L_{j},L_{i+k})\mu^{m}.
\end{eqnarray} Similarly, the value $a_2(L_i,L_j)$ depends only on the sum $i+j$. By again the skew-symmetry in \eqref{701}, $a_2(L_i,L_j)=0$ for $i,j\in\Z$. Following the procedure of Case \ref{caseee-1} we obtain  that $\phi_{\lambda}(L_{i},L_{j})=A(i+j)\lambda+B(i+j)\lambda^3$ (see also \cite{H}),    where $A$ is also a complex function .

\end{proof}

\end{document}